\documentclass[10pt]{article}
\usepackage[english]{babel}
\usepackage{amsfonts}
\usepackage[dvips]{graphicx}

\begin{document}

\title{\bf Numerical resolution of some BVP using Bernstein polynomials}

\date{October 2005}

\author{\bf Gianluca Argentini \\
\normalsize gianluca.argentini@riellogroup.com \\
\textit{Research \& Development Department}\\
\textit{Riello Burners}, 37048 San Pietro di Legnago (Verona), Italy}

\maketitle

\begin{abstract}
In this work we present a method, based on the use of Bernstein polynomials, for the numerical resolution of some boundary values problems. The computations have not need of particular approximations of derivatives, such as finite differences, or particular techniques, such as finite elements. Also, the method doesn't require the use of matrices, as in resolution of linear algebraic systems, nor the use of like-Newton algorithms, as in resolution of non linear sets of equations. An initial equation is resolved only once, then the method is based on iterated evaluations of appropriate polynomials.
\end{abstract}

\section{Basic concepts}

Let $B_{n,i}$ the $i$-th $n$-degree {\it Bernstein polynomial} (see \cite{bern}):

\begin{equation}\label{bernstein}
	B_{n,i}(t)= \left(_i^n\right)t^i(1-t)^{(n-i)}, \hspace{0.5cm} t\in[0,1]
\end{equation}

\noindent Let ${\bf P}_i=(p_i,q_i)\in{\mathbb R}^2$, with $0\leq i\leq n$. Then we can consider, from $[0,1]$ to ${\mathbb R}^2$, the {\it B$\acute{e}$zier curve} (see \cite{nurbs}) spawned by the array of points $({\bf P}_i)_{0\leq i\leq n}$:

\begin{equation}\label{bezier}
	bz(t)=\sum_{i=0}^n B_{n,i}(t){\bf P}_i
\end{equation}

\noindent Let $y=f(x)$ a real differentiable function define on an interval $[a,b]$ of ${\mathbb R}$, and let $(\psi(t),\phi(t))$, with $t\in[0,1]$, a possible parametric representation for $f$. Then, if $t\in [0,1]$, we have $x=\psi(t)$, $y(x)=f(\psi(t))=\phi(t)$ and, if $x_0=\psi(t_0)$, $\phi'(t_0)=f'(x_0)\psi'(t_0)$ (see \cite{diffg}). Therefore, if $\psi'(t_0)\neq 0$,

\begin{equation}\label{parametricFirstDerivation}
	\frac{df}{dx}(x_0)=\frac{\phi'(t_0)}{\psi'(t_0)}
\end{equation}

\noindent If $f\in C^k([a,b])$, than the B\'ezier curve

\begin{equation}
	bz(t)=\sum_{i=0}^n B_{n,i}(t)\left(\frac{i}{n},f\left(\frac{i}{n}\right)\right)
\end{equation}

\noindent and its derivatives until the $k$-th order are uniformely convergent to $f$ and its derivatives for $n\rightarrow +\infty$.\\
For our purposes, we'll transform a differential equation $F(y^{(k)}(x),x)=y(x)$, $x\in[a,b]$, from cartesian to parametric representation. As example, using the above notations and relations, the first degree ODE $y'(x)=x+y^2(x)$, supposed $\psi'(t)\neq 0$ in $(0,1)$, can be written in this form:

\begin{equation}
	\frac{\phi'(t)}{\psi'(t)} = \psi(t) + (\phi(t))^2
\end{equation}

\noindent Boundary values $y(a)=y_a$, $y(b)=y_b$ can be written in the form $\phi(0)=\phi_a=y_a$, $\phi(1)=\phi(b)=y_b$, where $a=\psi(0)$ and $b=\psi(1)$.\\
A result on the theorical resolution of an initial value problem for a particular ordinary differential equation is given in (\cite{approxB}).

\section{The algorithm of Bernstein polynomials}

Suppose we want a numerical solution for a BVP writable in the form

\begin{equation}\label{generalBVP}
	F(\phi^{(k)}(t), \psi(t)) = \phi(t),	\hspace{0.5cm} t\in [0,1], \hspace{0.5cm} \phi(0)=\phi_0, \phi(1)=\phi_1
\end{equation}

\noindent where $0\leq k\leq N$ for some $N\in{\mathbb N_+}$. Let $A$ the array of the {\it control points} $\{(a,ya),(b,y_b)\}$, and $dt$ the grid step for the parameter $t$.\\

{\bf 1}. As first part of the algorithm, let $(p,q)$ a pair of unknowns, and let $bz_1=(\psi_1(t),\phi_1(t))$ the B\'ezier curve spawned by the control points array $A=\{(a,y_a),(p,q),(b,y_b)\}$. Then we compute the values of $p$ and $q$ using the following two equations, obtained from the differential equation and from the boundary conditions respectively at $t=0$ and $t=1$:

\begin{eqnarray}
	\nonumber{ F(\phi_1^{(k)}(0), a) = y_a }\\
	\nonumber{ F(\phi_1^{(k)}(1), b) = y_b }
\end{eqnarray}

\noindent This system can have more then one real solution; in this case we choose the pair with the value of $p$ next to $0.5(b+a)$, the medium point of the cartesian interval $[a,b]$. We'll call the pair $(p,q)$ the {\it pivot} point. If the system have not real solutions, this first step has no output.\\
The pivot is useful to obtain a first raw graphic profile of the solution $y(x)$, since the line from $(a,y_a)$ to $(p,q)$ is tangent to the B\'ezier curve at $(a,y_a)$, and similarly for $(b,y_b)$ and $(p,q)$ (see \cite{nurbs}).\\

{\bf 2}. As second part of the algorithm, we open a loop starting from $m=1$. We evaluate the last B\'ezier curve $bz_m=(\psi(t),\phi(t))$, spawned by the actual control points array $A$, at the values $t_0 = (m \hspace{0.05cm} dt)$ and $t_1 = (1-m \hspace{0.05cm} dt)$ of the parameter $t$. Note that $\psi(t)$ and $\phi(t)$ are polynomials. In this way we obtain the two points $(\psi(t_0),\phi(t_0))$ and $(\psi(t_1),\phi(t_1))$. Then we can compute the {\it local errors}

\begin{equation}
	e_0 = F(\phi^{(k)}(t_0), \psi(t_0)) - \phi(t_0), \hspace{0.5cm} e_1 = F(\phi^{(k)}(t_1), \psi(t_1)) - \phi(t_1)
\end{equation}

\noindent These errors are an estimate about the approximation of the B\'ezier curve $bz_m$ as solution of the BVP ({\ref{generalBVP}). With the values $e_0$ and $e_1$, we can use a predefined criterion for deciding to exit from the loop or not. In the present work we use the following simple but effective criterion: let $s_{new}=0.5(e_1+e_2)$ and $s_{old}$ the corresponding value from the previous loop step. If $|s_{new}|\geq|s_{old}|$, then we exit from loop, else we insert the new points $(\psi(t_0),\phi(t_0)+e_0)$,$(\psi(t_1),\phi(t_1)+e_1)$ in the correct place of the control points array $A$, and perform a new step.\\
In any case, the loop finishes when a predefined maximum number of iterations is reached or when the product $m \hspace{0.05cm} dt$ is greater then the half of the interval $[0,1]$.\\

\noindent This technique can be considered a particular case of the {\it collocation method} (see \cite{polyApprox}) for the numerical solution of differential problems.

\section{A BVP with a linear differential equation}

As first application of the algorithm, we consider the following BVP:

\begin{equation}\label{bvp1}
	y'(x) = y(x) \hspace{0.5cm} \forall x\in(0,2), \hspace{0.5cm} y(0)=1, \hspace{0.5cm} y(2)=e^2
\end{equation}

\noindent which has the analytical solution $y=e^x$. According to (\ref{parametricFirstDerivation}), the problem has the following form when referred to the new variable $t\in[0,1]$:

\begin{equation}\label{bvp1parametrized}
	\frac{\phi'(t)}{\psi'(t)} = \phi(t) \hspace{0.5cm} \forall t\in(0,1), \hspace{0.5cm} \phi(0)=1, \hspace{0.5cm} \phi(1)=e^2
\end{equation}

\noindent The first step of the algorithm gives as pivot the point $(0.56344,1.68501)$. With a parametric step equal to $0.1$, the algorithm exits only after $3$ iterations, with a set of $9$ control points. Then the analitical approximated solution of the problem is a B\'ezier curve of $8$-th degree. In the following table we compare the numerical values between the approximated solution and the exact solution $e^{x(t)}$ for $t\in\{0.1i\}_{0\leq i\leq 10}$.\\

\begin{tabular}{l l l}
\textbf{t} & \textbf{approx} & \textbf{exact}\\
\scriptsize{0.0} & \scriptsize{1.0} & \scriptsize{1.0}\\
\scriptsize{0.1} & \scriptsize{1.17821} & \scriptsize{1.19321}\\
\scriptsize{0.2} & \scriptsize{1.38942} & \scriptsize{1.42075}\\
\scriptsize{0.3} & \scriptsize{1.75235} & \scriptsize{1.77812}\\
\scriptsize{0.4} & \scriptsize{2.35041} & \scriptsize{2.34651}\\
\scriptsize{0.5} & \scriptsize{3.18769} & \scriptsize{3.15902}\\
\scriptsize{0.6} & \scriptsize{4.17184} & \scriptsize{4.14995}\\
\scriptsize{0.7} & \scriptsize{5.14898} & \scriptsize{5.14912}\\
\scriptsize{0.8} & \scriptsize{5.98714} & \scriptsize{5.99204}\\
\scriptsize{0.9} & \scriptsize{6.67677} & \scriptsize{6.67179}\\
\scriptsize{1.0} & \scriptsize{7.38906} & \scriptsize{7.38906}\\
\end{tabular}\\
\\
\begin{figure}[ht]\label{grBVP1}
	\begin{center}
	\includegraphics[width=7cm]{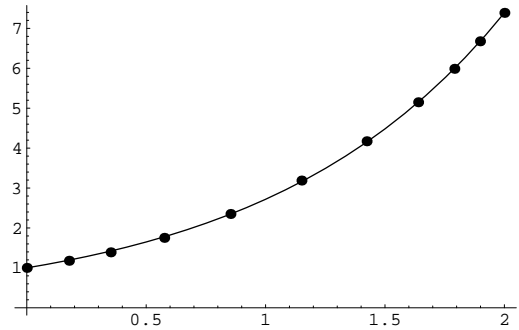}
	\caption{\small{\it The continuous line is the graphic of the exact solution, while dots are the cartesian pairs of the numerical solution as in previous table.}} 
	\end{center}
\end{figure}

\noindent In Fig.1 we compare the graphics of the two solutions.

\section{A BVP with a non linear differential equation}

In this section we consider the following BVP:

\begin{equation}\label{bvp2}
	-\frac{1}{4}\left(y'(x)\right)^2 + 4 = y(x) \hspace{0.5cm} \forall x\in(-1,3), \hspace{0.5cm} y(-1)=3, \hspace{0.5cm} y(3)=-5
\end{equation}

\noindent which has the analytical solution $y=-x^2+4$. In this case the solution has a local maximum in the interior of the domain. The problem has the following form when referred to the new variable $t\in[0,1]$:

\begin{equation}\label{bvp2parametrized}
	-\frac{1}{4}\left(\frac{\phi'(t)}{\psi'(t)}\right)^2 + 4 = \phi(t) \hspace{0.5cm} \forall t\in(0,1), \hspace{0.5cm} \phi(-1)=3, \hspace{0.5cm} \phi(3)=-5
\end{equation}

\noindent The first part of the algorithm gives as pivot the point $(1.0,7.0)$. With a parametric step equal to $0.1$, the algorithm exits only after $3$ iterations, with a set of $9$ control points. In the following table we compare the numerical values between the approximated solution and the exact solution $-x^2(t)+4$ for $t\in\{0.1i\}_{0\leq i\leq 10}$.\\

\begin{tabular}{l l l}
\textbf{t} & \textbf{approx} & \textbf{exact}\\
\scriptsize{-1.0} & \scriptsize{3.0} & \scriptsize{3.0}\\
\scriptsize{-0.721242} & \scriptsize{3.43889} & \scriptsize{3.47981}\\
\scriptsize{-0.439246} & \scriptsize{3.80217} & \scriptsize{3.80706}\\
\scriptsize{-0.0632671} & \scriptsize{4.02617} & \scriptsize{3.996}\\
\scriptsize{0.429831} & \scriptsize{3.8453} & \scriptsize{3.81525}\\
\scriptsize{1.0} & \scriptsize{3.02132} & \scriptsize{3.0}\\
\scriptsize{1.57017} & \scriptsize{1.56462} & \scriptsize{1.53457}\\
\scriptsize{2.06327} & \scriptsize{-0.226895} & \scriptsize{-0.257071}\\
\scriptsize{2.43925} & \scriptsize{-1.95482} & \scriptsize{-1.94992}\\
\scriptsize{2.72124} & \scriptsize{-3.44608} & \scriptsize{-3.40516}\\
\scriptsize{3.0} & \scriptsize{-5.0} & \scriptsize{-5.0}\\
\end{tabular}\\
\\
\begin{figure}[ht]\label{grBVP2}
	\begin{center}
	\includegraphics[width=7cm]{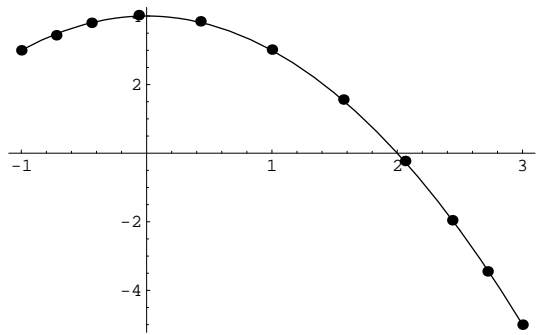}
	\caption{\small{\it The continuous line is the graphic of the exact solution, while dots are the cartesian pairs of the numerical solution as in previous table.}} 
	\end{center}
\end{figure}

\noindent In Fig.2 we compare the graphics of the two solutions.

\section{A more subtle BVP}

In this section we consider the following BVP:

\begin{equation}\label{bvp3}
	y'(x) = x + y^2(x) \hspace{0.5cm} \forall x\in(0,0.9), \hspace{0.5cm} y(0)=1, \hspace{0.5cm} y(0.9)=32.725
\end{equation}

\noindent which cannot be solved in terms of elementary functions (see \cite{diffEqBVP}). The built-in function {\ttfamily NDSolve} of the package {\it Mathematica 4.0} (see \cite{diffEqMath}) is not able to solve it directly, but we have numerically resolved the corresponding initial value problem for the condition $y(0)=1$, and then found the value for $y(0.9)$. Note that, in the interval $[0,0.9]$ considered, the derivative $y'$ is positive, hence the solution $y$ is increasing. With a parametric step equal to $0.15$, the algorithm of Bernstein polynomials exits after $2$ iterations, with a set of $7$ control points. In the following table we compare the numerical values between the approximated solution computed with this method ({\bf Bps}) and the approximated solution computed using {\ttfamily NDSolve}.\\

\begin{tabular}{l l l}
\textbf{x} & \textbf{Bps} & \textbf{NDSolve}\\
\scriptsize{0.0} & \scriptsize{1.0} & \scriptsize{1.0}\\
\scriptsize{0.13872} & \scriptsize{1.87331} & \scriptsize{1.17177}\\
\scriptsize{0.276777} & \scriptsize{2.77787} & \scriptsize{1.43186}\\
\scriptsize{0.418191} & \scriptsize{3.77639} & \scriptsize{1.85729}\\
\scriptsize{0.554867} & \scriptsize{4.80733} & \scriptsize{2.57996}\\
\scriptsize{0.674855} & \scriptsize{5.94924} & \scriptsize{3.84768}\\
\scriptsize{0.768333} & \scriptsize{7.55906} & \scriptsize{6.12075}\\
\scriptsize{0.831326} & \scriptsize{10.2845} & \scriptsize{10.0499}\\
\scriptsize{0.867159} & \scriptsize{14.9504} & \scriptsize{15.7566}\\
\scriptsize{0.885645} & \scriptsize{22.3191} & \scriptsize{22.256}\\
\scriptsize{0.9} & \scriptsize{32.725} & \scriptsize{32.725}\\
\end{tabular}\\
\\
\\
\noindent We note that locally the differences of values between the two methods can be large, but globally our method gives the same qualitative informations as the {\it Mathematica} internal function, in particular the growth in the interval $[0,0.9]$ and the presence of a value $x_0>0.9$ such that $\lim_{x\rightarrow x_0}y(x) = +\infty$. In Fig.3 we compare the graphics of the two solutions.

\begin{figure}[ht]\label{grBVP3}
	\begin{center}
	\includegraphics[width=7.5cm]{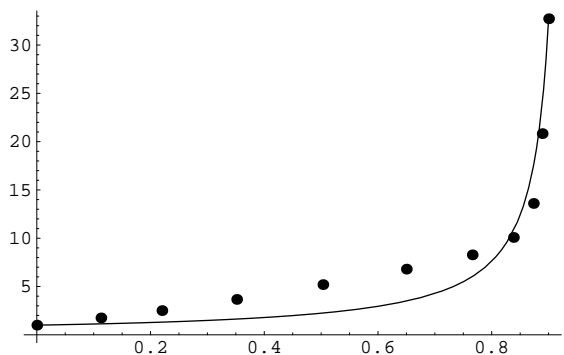}
	\caption{\small{\it The continuous line is the graphic obtained by the \textnormal{Mathematica} function {\textnormal {\ttfamily NDSolve}}, while dots are the cartesian pairs of the solution provided by the method of Bernstein polynomials.}} 
	\end{center}
\end{figure}

\section{Technical notes}

All the computations have been made using the software {\it Mathematica}, ver. 4.0, of Wolfram Research, on an Intel Xeon 3.2 GHz system. For an introduction to the use of this package in symbolic and numerical resolution of differential problems, see (\cite{diffEqMath}).\\
The basic Bernstein polynomials are build with the instructions\\

\noindent {\ttfamily Bernstein[i\_,n\_]:=Binomial[n,i]t\^{}i(1-t)\^{}(n-i);\\
Bernstein[n\_]:=Table[Bernstein[i,n],\{i,0,n\}]}\\

\noindent The symbolic derivative of Bernstein polynomials are computed with\\

\noindent {\ttfamily D[Bernstein[i,n],t]}\\

\noindent The B\'ezier curve are obtained using dot product:\\

\noindent {\ttfamily {\small Bezier[controlPoints\_]:=Bernstein[Length[controlPoints]-1].controlPoints}}\\

\noindent The computation of the pivot point is made using the built-in function {\ttfamily NSolve} for the resolution of a system of two equations. The evaluations of the B\'ezier curves at $t$-values are made using the transformation rules {\ttfamily /.} operator. \\

\noindent Electronic copies of the {\it notebooks} used for the numerical examples described in this work can be requested at the mail address {\it gianluca.argentini@gmail.com}.

\end{document}